\documentclass[leqno]{aupl}

\usepackage{
amsfonts,
latexsym,
amssymb,
enumerate,
verbatim,
mathrsfs,
graphicx,
centernot,
color,
rotate,}

\usepackage[all]{xy}

\newcommand{\labbel}[1]{\label{#1} [[{\bf #1}]]}  
\renewcommand{\labbel}{\label}

\newtheorem{theorem}{Theorem}[section]
\newtheorem{lemma}[theorem]{Lemma}

\newtheorem{corollary}[theorem]{Corollary}

\newtheorem*{claim*}{Claim}

\newtheorem*{theorem*}{Theorem}
\newtheorem*{proposition*}{Proposition}
\newtheorem*{corollary*}{Corollary}
\newtheorem*{lemma*}{Lemma}
\newtheorem*{scholion*}{Scholion}

\theoremstyle{definition}

\theoremstyle{remark}
\newtheorem{remark}[theorem]{Remark}

\newtheorem*{remark*}{Remark}
\newtheorem*{remarks*}{Remarks}

\newtheorem{conventions}[theorem]{Conventions}

\newtheorem*{observation*}{Observation}

 \allowdisplaybreaks[1]

\numberwithin{equation}{section}

\begin{document}

\title{Pairs of partial orders and the amalgamation property}

\author{Paolo Lipparini} 
\address{Dipartimento di $\overline{\text{A}}$malgamatematica\\Viale 
della  Ricerca
 Scientifica\\Universit\`a di Roma ``Tor Vergata'' 
\\I-00133 ROME ITALY
\\https://orcid.org/0000-0003-3747-6611}

\email{lipparin@axp.mat.uniroma2.it}

\urladdr{http://www.mat.uniroma2.it/\textasciitilde lipparin}

\subjclass{03C52; 06A75; 06F99}

\keywords{Strong amalgamation property; 
superamalgamation property; Fra\"\i ss\'e limit;
pairs of
binary relations; auxiliary relation; Urquhart doubly ordered sets; causal space}

\thanks{Work performed under the auspices of G.N.S.A.G.A. 
The author acknowledges the MIUR Department Project awarded to the
Department of Mathematics, University of Rome Tor Vergata, CUP
E83C18000100006.}

\begin{abstract}
We show that the theories of  partially ordered sets, 
lattices, semilattices, Boolean algebras, Heyting algebras
with a further coarser partial order, or a linearization, or an auxiliary
relation have the strong amalgamation property, 
Fra\"\i ss\'e limits and, in many cases, an $ \omega$-categorical
 model completion with quantifier elimination.
The same applies to Kronheimer and Penrose's causal spaces. 
On the other hand, 
Urquhart doubly ordered sets do not have the amalgamation property.

Our main tool is the superamalgamation property
(the strong amalgamation property is not enough), thus we provide
further arguments suggesting the usefulness of the superamalgamation property
also in pure model theory, not only in algebraic logic. 
\end{abstract}

\maketitle  

\section{Introduction} \labbel{in} 

The theories of linearly ordered sets
and of graphs are classical examples of theories 
in a  relational language and with the
 amalgamation property; they are basic examples
of theories allowing  Fra\"\i ss\'e limits for finite models.
Their Fra\"\i ss\'e limits are,
respectively, the ordered set of the rationals and
the random graph \cite{F,H}.
Partially ordered sets have the amalgamation property
and Fra\"\i ss\'e limits, as well
\cite{F,J}. 

More generally, consider the following properties of 
a binary relation $R$: 1.
$R$ is transitive;
2.
$R$ is reflexive;
3.
$R$ is symmetric;
4.
$R$ is antireflexive;
5.
$R$ is antisymmetric.
In \cite{duerel} we 
checked that the classic arguments show that,
  for every  $P \subseteq \{ 1,2,3,4,5 \} $,
the theory of  
a binary relation  satisfying the properties 
in $P$ has the  amalgamation property. 
Moreover, as reported in the following theorem,
we proved that the result extends
to pairs of comparable binary relations, possibly satisfying distinct sets
of properties.

\begin{theorem} \labbel{due} \cite[Theorem 3.1]{duerel} 
  For every pair  $P, Q \subseteq \{ 1,2,3,4,5 \} $,
the theory of  
a binary relation $R$ satisfying the properties 
from $P$, and of  a coarser
relation $S$ satisfying the properties 
from $Q$
 has the strong amalgamation property,
the superamalgamation property,
hence a model completion and a Fra\"\i ss\'e limit for finite
models. 
 \end{theorem}

Theorem \ref{due} does not generalize to triplets of
comparable relations \cite[Proposition 3.5]{duerel}; 
however, it does generalize to
arbitrary sets of \emph{transitive}  relations with any specified
family of comparability relations \cite{apu}.

Here we present  generalizations of Theorem \ref{due}
dealing with classes with the superamalgamation property
and with a further binary relation.
In particular, 
lattices, semilattices, 
Boolean algebras and Heyting algebras with a further coarser partial order,
or a linearization
have the amalgamation property.

We also consider some classical theories whose
 amalgamability status is not covered by the above results.
If $\mathbf S$ is a set with a partial order $\leq$ 
(henceforth, \emph{poset}, for short), an  \emph{auxiliary  relation}
 \cite[Definition I-1.11]{CLD}  is  a binary
antisymmetric and transitive relation $\ll$ 
which is finer than the order and satisfies
the following condition
\begin{equation}\labbel{aintr} \tag{A} 
w \leq x\ll y \leq z \text{ implies }  w \ll z.
 \end{equation}    
 The typical example 
of an auxiliary  relation  is the \emph{way below relation} 
\cite[I-1]{CLD}, which is an important 
tool in the theory of continuous lattices and domains,
which,  in turn, have many applications to the
 theory of computation,
 to the semantics of programming languages
and to many branches of  mathematics \cite{CLD}.
Auxiliary relations in ordered sets
with further structure arise
also in different contexts,
related to topological dualities, algebraic logic
and various generalizations of topology.
See \cite{CJ} for a survey.  

In Sections \ref{cau} and \ref{4}   we show that the theories of   
  posets (or lattices, semilattices, Boolean algebras, 
Heyting algebras) with an auxiliary
 relation have the strong amalgamation property and, in many cases, 
Fra\"\i ss\'e limits for finite structures, as well as model completions.
Similar results hold for the  theory of causal spaces,
introduced by 
Kronheimer and Penrose in \cite{KP}
in connection with foundational problems 
in general relativity.
In an equivalent formulation,
a \emph{causal space} is a poset 
with an antireflexive auxiliary relation.  
The formal similarities between the theory
of  posets with an auxiliary relation and
the theory of causal spaces
might have a deeper meaning;
see \cite{P} for a discussion. 
See Section \ref{cau} for more
details about the above notions.

More explicitly, our main results about the above-mentioned  notions 
are stated in the following theorem.

The order relation in an ordered structure
shall be usually denoted by $\leq$.
In particular, $\leq$ denotes the order naturally induced
by a join-semilattice, meet-semilattice, lattice or
Boolean structure.
In order
to have the joint embedding property, we assume that 
 Boolean algebras and Heyting algebras
are nontrivial (that is $0 \neq 1$) and that auxiliary relations in
Boolean algebras and Heyting algebras satisfy 
$0 \ll 0$ and  $1 \ll 1$.

\begin{theorem} \labbel{cp}
Suppose that $T^-$ is either the
theory of  partial orders, 
or  lattices, or  join semilattices, or  meet semilattices,
or Boolean algebras, or Heyting algebras.
In a language with an added binary relation $\ll$,
let $T$ be any extension of $T^-$
obtained  by adding one of the following axioms:
 \begin{enumerate}    
\item
$\ll$ is a reflexive order relation  coarser than $\leq$, or
\item
$\ll$ is a linearization of $\leq$, or
\item
$\ll$ is an auxiliary relation.
  \end{enumerate}
Alternatively, let $T$ be   the theory of causal sets.

Then $T$  has the strong amalgamation property,
more generally, the superamalgamation property 
with respect to $\leq$. 

The class of finite 
models of $T$ has
 a Fra\"\i ss\'e limit $\mathbf M$. 
Except possibly for lattices and Heyting algebras, 
the first-order theory of $\mathbf M$ is 
$ \omega$-categorical, has quantifier 
elimination and is  the model completion
of $T$.  
 \end{theorem} 

Parts (1) and (2) of Theorem \ref{cp}
will be proved at the end of Section \ref{2}.  
In the case of partial orders
Theorem \ref{cp}  (3)
is a special case of
Theorem \ref{gen} 
which will be proved in Section \ref{cau}.
The same applies  to the theory  of causal sets. 
The remaining cases follow from Theorem \ref{genn} in Section \ref{4}.
Details  are given at the end of Section \ref{4}.

What is relevant in the proof of Theorem \ref{cp}
is that lattices, semilattices, etc. have the 
superamalgamation property: generally, the result applies
to any such theory. See Theorem \ref{genn} below. 
The superamalgamation property is necessary in the arguments:
the strong amalgamation property alone is not enough.
See Remarks \ref{rmk} and \ref{rmkk}.  
So far, the main applications of the superamalgamation property 
have been found in algebraic logic, e.~g., \cite{GM,KH,Ma}.
Our results show that the superamalgamation property has also some 
model-theoretical interest. 
A similar situation occurred in \cite{sapimpap},
where the reader might find more comments. 

Considering another theory
with two binary relations, doubly ordered sets
have been introduced by 
Urquhart \cite{U} 
in connection with the representation theory of lattices.
In Section \ref{Usec} we show that 
 Urquhart doubly ordered sets do not have the amalgamation property;
more generally, the amalgamation property fails for pairs of transitive
relations satisfying the   Urquhart condition.

\smallskip 

We now recall the basic notions. 
\emph{Models} are intended in the classical 
model-theoretical sense \cite{H}. 
An \emph{ordered structure} is a model with 
a partial order $\leq$ and, possibly, further relations, functions
and constants. 

A class $\mathcal K$ of models  
of the same type has the \emph{amalgamation property}
(AP)  if, whenever 
$\mathbf A, \mathbf B, \mathbf C \in \mathcal K$,
 $ \iota_1  \colon \mathbf C \to \mathbf A$
and 
 $ \kappa _1  \colon \mathbf C \to \mathbf B$
are embeddings, then there are a model
$\mathbf D \in \mathcal K$ and  embeddings
$ \iota  \colon \mathbf A \to \mathbf D$
and 
 $ \kappa   \colon \mathbf B \to \mathbf D$
such that 
$  \iota_1  \circ \iota
=
\kappa_1   \circ \kappa $. 
If, in addition,
$\mathbf  D$,
$ \iota  $
and 
 $ \kappa   $ can be always chosen in such a way that 
$\iota(A) \cap \kappa  (B) = (\iota_1 \circ \iota)(C)$,
then $\mathcal K$ is said to have the \emph{strong amalgamation property}
(SAP).
See \cite{KMPT} for a comprehensive survey of classical results
about AP.

Suppose further that $\mathcal K$ is a class of ordered structures with order 
$\leq$. Then $\mathcal K$ has the \emph{superamalgamation property
with respect to $\leq$,}
or simply the \emph{superamalgamation property}, when $\leq$ is understood,
if, in addition, for every $a \in A \setminus \iota _1(C) $ and 
$b \in B \setminus \kappa  _1(C)$,
if $ \iota(a) \leq _{ \mathbf  D} \kappa (b)$,
then there is $c \in C$ such that 
$a \leq _{ \mathbf  A} \iota_1 (c)$,
$ \kappa_1(c) \leq _{ \mathbf  B} b$, 
and also the corresponding conclusion holds when 
$ \iota(a) \geq _{ \mathbf  D} \kappa (b)$. The superamalgamation property 
has found significant applications
in algebraic logic \cite{GM,KH}.
Notice that, since
$\leq$  is transitive,
the superamalgamation property  determines
$\leq$  uniquely
on $\iota_1 (\mathbf A) \cup \kappa _1 (\mathbf  B)$. 

We shall use the above definition also
for an arbitrary binary relation $R$  in place of the order
$\leq$. The generalization has some use, as well \cite{Ma}.

If $\mathcal H$ is another class of 
models of the same type and, for every $\mathbf A, \mathbf  B,
 \mathbf  C $ in $ \mathcal K$ and $ \iota_1, \kappa _1$ as above,  
there exist some $\mathbf D $ in $  \mathcal H$---not
 necessarily in $  \mathcal K$---and embeddings
 $ \iota, \kappa $ as above,
 we say that 
\emph{$\mathcal K$ has the (strong, super) amalgamation property in $\mathcal H$}.

\section{Ordered structures with the superamalgamation property} \labbel{2} 

In this section we give a proof of 
Parts (1) - (2) of Theorem \ref{cp}.
The section is a good introduction to the methods
which will be used in the subsequent sections.

\begin{lemma} \labbel{lemsolo}
Suppose that $(D, {\leq_D},{\ll_D})$
is a set with two partial orders, and $\ll_D$ is
coarser than  $\leq_D$. 
If $(E, {\leq_E})$ is a poset extending
 $(D, {\leq_D)}$, 
 then 
there is an order $\ll_E$ on $E$ such that
$\ll_E$ is coarser than $\leq_E$ and 
 $(E, {\leq_E},{\ll_E} )$ extends $(D, {\leq_D},{\ll_D})$.
 \end{lemma} 

\begin{proof}
Let 
$e \ll_ { E} f$
on $ E$ if  
either $e \leq_ { E} f$, or
 there are $h,k \in D$
such that  $e \leq_ { E} h \ll_ { D} k \leq_ { E} f$.
By construction, $\ll_ { E}$ is coarser than $\leq_ { E}$.  
Moreover, $\ll_ { E}$ extends $\ll_ { D}$.
Indeed, if $e,f \in D$ and  $e \ll_ { D} f$,
then $e \ll_ {E} f$ by definition and since $\leq_ { E} $
is reflexive. In the other direction, 
if $e,f \in D$ and  $e \ll_ { E} f$ is given by 
$e \leq_ { E} f$, then $e \leq_ { D} f$,
since $e,f \in D$ and $\leq_ { E} $ extends $ \leq_ { D} $,
thus $e \ll_ { D} f$, since by assumption
$\ll_ { D} $ is coarser than $ \leq_ { D} $. 
On the other hand, if $e,f \in D$ and  $e \ll_ { E} f$ is given by 
$e \leq_ { E} h \ll_ { D} k \leq_ { E} f$, then
$e \leq_ { D} h \ll_ { D} k \leq_ { D} f$,
since $e,h,k,f \in D$. Thus
$e \ll_ { D} h \ll_ { D} k \ll_ { D} f$,
since
$\ll_ { D} $ is coarser than $ \leq_ { D} $,
then $e \ll_ { D} f$ by transitivity of 
$\ll_ { D}$.

We now check that $\ll_ { E} $
is transitive.
If $e \ll_ { E} f \ll_ { E} g$ is witnessed by
\begin{equation}\labbel{wit}    
e \leq_ { E} h \ll_ { D} k \leq_ { E} f \leq_ { E} p \ll_ { D} q \leq_ { E} g,
   \end{equation} 
then $k \leq_ { E} p$, by transitivity of $\leq_ { E}$. 
Since $k,p \in D$, then  $k \leq_ { D} p$,
hence $k \ll_ { D} p$, since $ \ll_ { D}$ 
is coarser than $\leq_ { D} $. Hence
$h \ll_ { D} k \ll_ { D} p \ll_ { D} q$,
thus $h \ll_ { D} q$,
by transitivity of $\ll_ { D} $.
Then \eqref{wit} reads  
$e \leq_ { E} h \ll_ { D}  q \leq_ { E} g$,
which means $e \ll_ { E} g$ by definition.
The other cases are immediate from
transitivity of $\leq_ { E}$.

It remains to check that 
 $\ll_ { E}$ is antisymmetric.
Suppose that
$e \ll_ { E} f \ll_ { E} e$ is witnessed by
\eqref{wit} with $g=e$.   
We have proved above that 
$h \ll_ { D} k \ll_ { D} p \ll_ { D} q$.
Moreover, $q \leq_ { E} g=e \leq_ { E} h$,
hence  $q \leq_ { E} h$, $q \leq_ { D} h$
and $q \ll_ { D} h$, since $q,h \in D$
and $ \ll_ { D} $ is coarser than $\leq_ { D} $.
Thus $h \ll_ { D} k \ll_ { D} p \ll_ { D} q \ll_ { D} h$,
hence $h=k=p=q$, by transitivity and antisymmetry
of $\ll_ { D}$. Finally,   
from $h=k $, $ p=q$ and \eqref{wit} with $g=e$ we get
$e \leq_ { E} h \leq_ { E} f \leq_ { E} p  \leq_ { E} e$,
hence   $e=f$, by   by transitivity and antisymmetry
of $\leq_ { E}$. The other cases are much simpler.
For example, if $e \leq_ { E} h \ll_ { D} k \leq_ { E} f \leq_ { E} e$,
then  $k \leq_ { E} f \leq_ { E} e \leq_ { E} h $, hence
$k \leq_ { E} h $,  
$k \leq_ { D} h$,  $k \ll_ { D} h$, thus $k=h$ by
antisymmetry  of $\ll_ { D}$. Hence 
$e \leq_ { E} h   \leq_ { E} f \leq_ { E} e$, thus $e=f$. 
 \end{proof}

The next theorem collects some classical results,
sometimes  not explicitly stated.
The theorem follows from
some proofs in \cite{GM,J,Ko}. 
See the proof of \cite[Theorem 2.4]{sapimpap} 
for full details.

\begin{theorem} \labbel{josup}
The classes of partially ordered sets,  meet semilattices,  join semilattices, 
lattices,  Boolean algebras and Heyting algebras
 have the superamalgamation property.
The same applies to the classes of finite such structures.
 \end{theorem}

\begin{remark} \labbel{un}
Assume for simplicity that the structures to be amalgamated
are such that $A \cap B =C$ and that 
the embeddings $\iota_1$ and $\kappa_1$ 
are inclusions. In the proof of
Theorem \ref{due}, in each case,
the amalgamating model has been constructed over
$D=A \cup B$. Independently from the actual proof,
this fact follows anyway from the statement, since Theorem \ref{due} 
deals with relational languages, and if some universal relational theory has
the strong (super) amalgamation property, then it has (super)amalgamation
``over union'': just consider an appropriate substructure. 
A similar remark  applies to
Theorem \ref{gen} below.
 
The possibility of having the amalgamating model  over
$D=A \cup B$ sometimes proves  to be useful \cite{apu}.
As additional examples, the present remark will be also used in
the next proof, in
the proof of Theorem \ref{genn} 
and in Remark \ref{add}(a) 
below. 
 \end{remark}

\begin{proof}[Proof of Theorem \ref{cp}, Part (1), Superamalgamation]
Without loss of generality, assume that 
$\mathbf A$, $\mathbf  B$, $\mathbf  C$ 
are members of $T$ such that  
$\mathbf  C \subseteq \mathbf A, \mathbf  B$,   
and $A \cap B =C$. 
Let  $\mathbf A^-$, $\mathbf  B^-$ and $\mathbf  C^-$ 
be the reducts of  $\mathbf A$, $\mathbf  B$
and $\mathbf  C$ to the language of $T^-$.
By Theorem \ref{josup}, there is a model $\mathbf E^-$
 superamalgamating 
  $\mathbf A^-$ and $\mathbf  B^-$ 
over $\mathbf  C^-$ in the appropriate class
(posets, semilattices, etc.). Without loss of generality,
we may assume that 
$\mathbf E^-$ extends both $\mathbf A^-$ and $\mathbf  B^-$.
Let $D=A \cup B$. 
By Theorem \ref{due} and Remark \ref{un},   
the $\{ {\leq},{\ll}\}$-reducts of $\mathbf A$ and $\mathbf  B$ 
can be superamalgamated over 
the $\{ {\leq},{\ll}\}$-reduct of
$\mathbf  C$ by an amalgamating structure $(D, {\leq},{\ll})$
over $D$. In particular, say, $(A, {\ll})$
embeds in  $(D, {\ll})$. 
 
Since the superamalgamation property 
determines the order on $A \cup B$ uniquely,
then  $(D, {\leq})$ is a substructure of $(E, {\leq})$.
By Lemma \ref{lemsolo}, there is an order $\ll_E$ on $E$ 
such that $(D, {\leq},{\ll})$ embeds in $(E, {\leq},{\ll_E})$
and 
$\ll_E$ is coarser than $\leq_E$.
If we expand $\mathbf E^-$ by adding $\ll_E$,
we get a model $\mathbf E$ of $T$.
Since $(A, {\ll})$
embeds in  $(D, {\ll})$ and 
$(D, {\ll})$ embeds in $(E,{\ll_E})$,
then $(A, {\ll})$ embeds in $(E,{\ll_E})$,
and the same for $(B, {\ll})$.
Since $\mathbf E^-$ extends both $\mathbf A^-$ and $\mathbf  B^-$,
this takes care of  the language of $ T^-$.
Thus $\mathbf E$  superamalgamates 
  $\mathbf A$ and $\mathbf  B$ 
over $\mathbf  C$.

\emph{(Part 2, superamalgamation)} If $\ll$ is a linearization
of $\leq$ in $\mathbf A$, $\mathbf  B$, $\mathbf  C$, then,
in particular, $\ll$ is an order coarser than $\leq$.
We can thus apply (1) in order to get a model $\mathbf E$  amalgamating
$\mathbf A$ and $\mathbf  B$ over $\mathbf  C$,
where $\ll_E$ in $\mathbf E$  is a
(not necessarily linear)  order coarser than $\leq_E$.
 Let $ \mathbf F$ be obtained from $\mathbf E$ 
  by replacing $\ll_E$ with a linearization $\ll_F$ of $\ll_E$.
The only thing to check is that $ \mathbf F$ extends $\mathbf A$ and $\mathbf  B$.
If $a, b \in A$, $a \neq b$  and $a \ll_F b$, then it is not the case that
$b \ll_A a$, since $\ll_E$ extends $\ll_A$, and $\ll_F$
is antisymmetric and coarser than $\ll_E$. Since $\ll_A$ is a linear order,
then $a \ll_A b$. The case $a \neq b \in B$ is treated in the same way.
Thus $ \mathbf F$ amalgamates $\mathbf A$ and $\mathbf  B$
over $\mathbf  C$.

\emph{(Parts (1) - (2), Fra\"\i ss\'e limits, quantifier elimination and model-completion)}
The remaining statements in the theorem  follow from 
the just proved amalgamation property using
 standard model-theoretical arguments.

In the cases of partial orders, semilattices and
lattices,
the joint embedding property
follows from the amalgamation property, since
we are allowed to consider $\mathbf  C$ as an empty structure.
In the cases of Boolean  and Heyting algebras,
by the assumptions stated right before Theorem \ref{cp},
 there is a unique $ \emptyset $-generated structure,
 which can be taken as $\mathbf  C$, thus the amalgamation property 
provides the joint embedding property in such cases, as well.
Notice that if $\mathbf A$, $\mathbf  B$ and $\mathbf  C$ 
are finite, then, in view of the last statement in 
Theorem \ref{josup}, the above proof provides a finite amalgamating
model $\mathbf  D$ in each case. Hence the class of finite models
of $T$ has the amalgamation property and the joint 
embedding property. 

Then use  \cite[Chapter 7]{H},
specifically,
Theorems 7.1.2 and 7.4.1 therein. 
As far as the statement about the model completion is concerned,
see \cite[Fact 2.1(3)]{KS};
the proof here does not apply to lattices and Heyting algebras,
since such theories are not locally finite.
\end{proof}

\section{Posets with an auxiliary  relation; causal spaces} \labbel{cau}

If $\mathbf S= (S, {\leq})$ is a poset, an \emph{auxiliary order},
or an \emph{auxiliary  relation} \cite[Definition I-1.11]{CLD} 
on $\mathbf S$ is  a binary
 relation $\ll$  finer than the order and satisfying
\begin{equation}\labbel{a} \tag{A} 
w \leq x\ll y \leq z \text{ implies }  w \ll z,
 \end{equation}    
for all $w, x,  y, z \in S$. Since $\ll$ is finer than $\leq$
and $\leq$ is antisymmetric, 
then $\ll$ is antisymmetric, as well.
Notice that we neither require $\ll$
to be reflexive, nor to be antireflexive.
On the other hand, $\leq$ is assumed to be reflexive.
It follows that \eqref{a} is equivalent to the conjunction of 
the following two conditions.
\begin{align} \labbel{a1}    \tag{A1}
&w \leq x\ll y \text{ implies }  w \ll y, 
\\
\labbel{a2}    \tag{A2}
& x\ll y \leq z \text{ implies }  x \ll z.
 \end{align} 
Since $\ll$ is finer than $\leq$,
it follows that $\ll$  is transitive.

As we mentioned, we do not require $\ll$
 to be antireflexive in the definition of an auxiliary relation.
If we further require that $\ll$ is antireflexive, then the structure
 $(S, {\leq}, {\ll})$ is  
a \emph{causal space} in the terminology from \cite{KP}.
If this is the case, then
$\ll$ is a strict  partial order.
In \cite{KP} another relation $\rightarrow $   
is considered, but it can be defined in terms of 
$\leq$ and $\ll$ by a universal sentence, hence
$\rightarrow $ does not modify the notion of an embedding
(it  modifies the notion of a homomorphism,
however). 

In all the above definitions we might assume that 
$\mathbf S $
is a lattice
$ (S, {\wedge}, {\vee})$,
or   a join semilattice
$ (S,  {\vee})$,
or a meet semilattice
$ (S, {\wedge})$.
Algebraically, there is no difference between
join and meet semilattices; however,
in the former case $a \leq b$
is defined by  $ a \vee b = b$,
while in the latter case $a \leq b$
is defined by  $ a \wedge b = a$.
In the case of lattices the above definitions of $\leq$ 
are equivalent.
In each case, the definition of an auxiliary relation
is given with reference to the order $\leq$ as introduced above.

The way back relation---the typical example
of an auxiliary relation---is frequently considered in 
lattice-ordered structures \cite{CLD}; on the other hand the lattice-ordered
version of causal sets is possibly  deprived 
of physical sense. We shall also consider Boolean algebras and Heyting algebras
endowed with an auxiliary relation. 

Condition \eqref{a} has been considered in models 
with further structure, 
frequently under the name \emph{subordination}, e.~g., clause (S4)
in either \cite[Definition 2.1]{BBSV}, or  \cite[Definition 9]{C}, or \cite[Definition 1]{CJ}.
See the quoted papers for credits 
to original sources and further references. 
 
Let us mention that Clause \eqref{a}
arises also from a very general situation.
Suppose that $\mathbf S$ is a poset,
$K$ is a unary operation on $S$ and set
$x \ll y$ if $Kx \leq y$.
Condition \eqref{a2} follows just from
transitivity of $\leq$. 
If $K$ is extensive, that is,
$x \leq Kx$ holds for every $x \in S$,
then   $x \ll y$ is finer than $\leq $.
 If furthermore $K$ is isotone, 
then \eqref{a1} holds.
In conclusion, if
$K$ is isotone and extensive, then 
 $(S, {\leq}, {\ll})$ is a poset with an
auxiliary relation.
 Dually, if $I$ is an isotone and 
contractive operation, then 
we get an auxiliary relation by setting
$x \ll y$ if $x \leq Iy$.
The above remarks are a variation on
known ideas, 
e.~g., \cite[Remark 2]{CJ}.

We now  embark on the proof of a 
generalization of Theorem \ref{cp}.
The proof will occupy the next two sections.
The present section essentially deals 
with posets with a further binary relation,
which might be a comparable order, an auxiliary relation or, 
possibly, a relation satisfying a set of distinct properties.
In the next section we shall first prove an extension Lemma \ref{lemsolo},
according to which an additional binary relation on a poset $\mathbf P$ 
can be lifted to a poset extending $\mathbf P$.
The superamalgamation property for all the theories considered
here then allows us to prove Theorem \ref{cp}
in a slightly more general form in which the additional 
relation $\ll$ is allowed to satisfy a range of possibilities.

\begin{conventions} \labbel{conv}    
We shall  consider the following properties of 
a binary relation $R$.
  \begin{enumerate}   
 \item 
[2.]
$R$ is reflexive;
\item
[4.]
$R$ is antireflexive;
\item
[5.]
$R$ is antisymmetric. 
  \end{enumerate} 
(the numbering is intended to be consistent 
with \cite{duerel,apu}; here all the relations are always
1.\ transitive and never assumed to be 3.\ symmetric.)
We will consider two binary relations
$\ll$ and $\leq$.
By F and C, respectively, we will mean the conditions that 
  $\ll$ is finer, coarser, respectively, than $\leq$.
\end{conventions}

\begin{theorem} \labbel{gen}
Suppose that $P, Q \subseteq \{ 2,4,5 \} $
and $N \subseteq \{F, C, A1, A2 \} $. 
Suppose further that either
  \begin{enumerate}[(a)]   
 \item 
 $\{ A1, A2 \}  \subseteq N$, or 
\item
$A1 \notin N$, $A2 \notin N$, or
\item   
both
(c1) $F \in N$ and
(c2) either  $5 \in P$, or $4 \in Q$,  or $5 \notin Q$. 
  \end{enumerate} 

Let $T= T _{P,Q,N} $ be the theory asserting that 
$\ll$ and $\leq$ are transitive binary relations, that
$\leq$ satisfies the properties in $P$, $\ll$
satisfies the properties in $Q$ and that the properties
in $N$ are satisfied.

Then $T$ has the superamalgamation property
with respect to $\leq$.
The class of finite 
models of $T$ has
 a Fra\"\i ss\'e limit $\mathbf M$.
The first-order theory of $\mathbf M$ is 
$ \omega$-categorical, has quantifier elimination and is  the model completion
of $T$.  
 \end{theorem}

Notice that, for certain combinations of
$P$, $Q$ and $N$, we get  a trivial conclusion, or a
trivial class of structures.
For example, if $F, C \in N$, then necessarily
$\ll$ is equal to $\leq$, hence there is nothing to prove
(provided the result is known for a single relation). 
If $2 \in Q$  and either  
$A1 \in N$ or $A2 \in N$, then necessarily
$\ll$ is coarser than $\leq$, thus if also
$F \in N$, then again $\ll$ is equal to $\leq$.
On the other hand, we cannot have both
$2$ and $4$ in $P$ (or in $Q$).
Similarly, we cannot have $2 \in P$,
$4 \in Q$ and $C \in N$. 
In order for the statement of Theorem \ref{gen} 
to be formally true in the last two cases,
 we 
allow empty structures, or else
the reader  should discard such cases.

 \begin{proof}
The cases when $A1 \notin N$ and $A2 \notin N$
are covered by \cite[Theorem 3.1]{duerel},
reported here as Theorem \ref{due}.
Formally, Theorem \ref{due} does not deal with the case
$F,C \notin N$, but in this case there are no connections
between $\leq$ and $\ll$, hence we can work
independently and the case of a single relation \cite[Proposition 2.1]{duerel}
is enough, applied twice.  

We now divide the proof into cases.

(I) $A1 \notin N$ and $A2 \in N$.
Hence, by assumption, $F \in N$,
that is,
$\ll$ is finer than $\leq$.  Moreover, (c2) holds.

Without loss of generality,
assume that the models to be amalgamated are 
$\mathbf  C \subseteq \mathbf A, \mathbf  B$   
with $A \cap B =C$. 

Define $\leq$ over  $D=A \cup B$ as
 the union of $\leq_ {\mathbf  A}$, $\leq_ {\mathbf  B}$,
${\leq_ {\mathbf  A}} \circ {\leq_ {\mathbf  B}}$
and ${\leq_ {\mathbf  B}} \circ {\leq_ {\mathbf  A}}$.
By the proof of  \cite[Proposition 2.1 (A), case b]{duerel}, $(D, \leq)$
superamalgamates the $\leq$-reducts of $\mathbf A$ and $\mathbf  B$
over $\mathbf  C$. This is a revisitation of \cite[Lemma 3.3]{J}.

Define $\ll$ over  $D=A \cup B$ as
 the union of $\ll_ {\mathbf  A}$, $\ll_ {\mathbf  B}$,
${\ll_ {\mathbf  A}} \circ {\leq_ {\mathbf  B}}$ and
${\ll_ {\mathbf  B}} \circ {\leq_ {\mathbf  A}}$,
and consider the model $\mathbf  D = (D, \leq,\ll)$.
In order to lighten the notation,
unsubscripted relations are always meant as interpreted in $\mathbf  D$.

We first show that, say, $\mathbf A$ is actually a substructure of 
$\mathbf  D$. The $\leq$-part follows from the mentioned proof
by \cite{J,duerel}. If $a,a' \in A$ and 
$ a \ll_ {\mathbf  A} a'$, then $ a \ll a'$, by the definition of 
$\ll$.   
In the other direction, suppose that $a,a' \in A$ and
$a \ll a'$ holds in $\mathbf  D$ because, say, $(a,a')$ belongs to 
${\ll_ {\mathbf  A}} \circ {\leq_ {\mathbf  B}}$,
that is,     
\begin{equation}\labbel{eq}     
 a \ll_ {\mathbf  A} c \leq_ {\mathbf  B} a',
  \end{equation}
for some $c \in D$.
Then $c \in A$, because of the $\ll_ {\mathbf  A}$ relation,
and $c \in B$, because of the $\leq_ {\mathbf  B}$ relation,
thus $c \in C$, since $C=A \cap B$.
Also, by the   $\leq_ {\mathbf  B}$ relation, 
$a' \in B$, thus  $a' \in C$.
Since both $c$ and  $a'$ belong to $C$,
then   $c \leq_ {\mathbf  C} a'$, since $\mathbf  B$ extends $\mathbf  C$.
But also $\mathbf A$ extends $\mathbf  C$, thus
$c \leq_ {\mathbf  A} a'$. Then \eqref{eq} reads 
$ a \ll_ {\mathbf  A} c \leq_ {\mathbf  A} a'$.
By \eqref{a2}, which holds in $\mathbf A$ by assumption,
we have   $ a \ll_ {\mathbf  A} a'$. The other cases are similar
or simpler. Similarly,  $\mathbf B$ is  a substructure of 
$\mathbf  D$.

We now show that  \eqref{a2} holds in  $\mathbf  D$.
The proof can be given by analyzing
all possible cases, but we will provide  a more uniform argument.
Suppose that $x \ll y$, $y\leq z$ and, say, $x  \in A$.
In any case, $x \ll y$ is witnessed by a relation of the form
$ x \ll_ {\mathbf  A} w$, for some $w \in D$, possibly 
followed by a $ \leq_ {\mathbf  A} $
or a $ \leq_ {\mathbf  B} $ relation added on the right.
For example, if  $x \ll y$ is given by $ x \ll_ {\mathbf  B} c \leq_ {\mathbf  A} w$,
for $c \in C$, then $x \in A \cap B=C$, hence also  
$ x \ll_ {\mathbf  A} c$, since $\mathbf  C$ embeds both in $\mathbf A$
and in $\mathbf  B$.
Moreover, $y\leq z$ is witnessed by $ \leq_ {\mathbf  A} $
and $ \leq_ {\mathbf  B} $ relations only. 
Thus if $x \ll y\leq z$, then $x$ and $z$ can be connected by a chain
of relations starting with  $ \ll_ {\mathbf  A} $ and with
 possibly further $ \leq_ {\mathbf  A} $s
and $ \leq_ {\mathbf  B} $s on the right. 
Now we show how to reduce the length of such 
chains of relations.
If
$ x \ll_ {\mathbf  A} w \leq_ {\mathbf  A} v$, 
then $ x \ll_ {\mathbf  A} v$ by \eqref{a2} in $\mathbf A$,
thus we may assume that no $ \leq_ {\mathbf  A} $  
appears immediately after  $\ll_ {\mathbf  A}$. 
If $ x \ll_ {\mathbf  A} w \leq_ {\mathbf  B} v \leq_ {\mathbf  B} u$,
then $ x \ll_ {\mathbf  A} w \leq_ {\mathbf  B} u$ by transitivity
of $\leq_ {\mathbf  B}$ in $\mathbf  B$. If
$ x \ll_ {\mathbf  A} w \leq_ {\mathbf  B} v \leq_ {\mathbf  A} u$,
then $w,v \in A \cap B$, hence, by an argument usual by now, 
$w \leq_ {\mathbf  A} v$.
The above relation then reads 
$ x \ll_ {\mathbf  A} w \leq_ {\mathbf  A} v \leq_ {\mathbf  A} u$,
hence 
$ x \ll_ {\mathbf  A} u$ by \eqref{a2} in $\mathbf A$. 
Iterating the above reductions, we get that if 
$x \in A$ and
$x \ll y\leq z$, then either 
$ x \ll_ {\mathbf  A} z$, or 
$ x \ll_ {\mathbf  A} w \leq_ {\mathbf  B} z$.
The case $x \in B$ is symmetrical. 
In each case 
$x \ll z$,
by the definition of $\ll $.

Since $\ll$  is finer than $\leq$
in both $\mathbf A$ and $\mathbf  B$, then 
 $\ll$  is finer than  $\leq$ in $\mathbf  D$.
Transitivity of $\ll$ then follows from the already proved
condition \eqref{a2}. 

If $\ll$ is antireflexive in $\mathbf A$ and $\mathbf  B$,
then  $\ll$ is antireflexive in $\mathbf  D$, since
$D= A \cup B$ and both $\mathbf A$ and $\mathbf  B$
embed in $\mathbf  D$, as we have already showed.
For the same reason, if $\ll$ is reflexive in $\mathbf A$ and $\mathbf  B$,
then  $\ll$ is reflexive in $\mathbf  D$,

Assume that $5 \in Q$,
that is, $\ll$ is supposed to be 
antisymmetric; in particular
$\ll$ is antisymmetric in $\mathbf A$ and $\mathbf  B$.
By the assumption (c2),
either   $4 \in Q$, that is, 
$\ll$ is antireflexive, or $5 \in P$, that is,
$\leq$ is  antisymmetric.
In the former case  $\ll$ is necessarily antisymmetric,
since  it is both
transitive and antireflexive.
In the latter case,
if $x \ll y $, $ y \ll x $,
then $x \leq y $, $ y \leq x $,
since $\ll$ is finer than $\leq$.
Thus $x=y$, by antisymmetry of $\leq$.

(II) $A1 \in N$ and $A2 \notin N$.
This case is symmetrical to (I).

(III) $A1 \in N$ and $A2 \in N$.
The proof is similar to the proof in (I);
however, the definition of $\ll$ in $\mathbf  D$ 
should be adapted to this case, and we need to show
how to get by without    
the assumptions (c1) and (c2).

Define $D$ and $\leq$ as in (I) above.
In the present  case define $\ll$ over  $D=A \cup B$ as
 the union of $\ll_ {\mathbf  A}$, $\ll_ {\mathbf  B}$,
${\ll_ {\mathbf  A}} \circ {\ll_ {\mathbf  B}}$,
${\ll_ {\mathbf  B}} \circ {\ll_ {\mathbf  A}}$,
${\ll_ {\mathbf  A}} \circ {\leq_ {\mathbf  B}}$,
${\leq_ {\mathbf  A}} \circ {\ll_ {\mathbf  B}}$,
${\ll_ {\mathbf  B}} \circ {\leq_ {\mathbf  A}}$ and
${\leq_ {\mathbf  B}} \circ {\ll_ {\mathbf  A}}$.

The proof that $\mathbf  D$ extends both $\mathbf A$
and $\mathbf  B$ is not significantly different from (I).
If, say, $a, a' \in A$ and $a \ll a'$ in $\mathbf  D$, then
this is witnessed either by  
$a \ll_ {\mathbf  A} a'$ (and there is nothing to prove),
or by 
$a \ll_ {\mathbf  B} a'$ (and then $a,a' \in C$, hence  $a \ll_ {\mathbf  A} a'$,
again),
or by  
relations of the form
$ a \trianglelefteq _ {\mathbf  A} c \blacktriangleleft_ {\mathbf  B} a'$
or
$ a \trianglelefteq _ {\mathbf  B} c \blacktriangleleft_ {\mathbf  A} a'$,
for appropriate $\trianglelefteq ,  \blacktriangleleft$.
 In any case, $c \in C$ and, say, in the former case, 
$a' \in B$, hence $a' \in C$, thus the relations can be evaluated in $\mathbf A$,
providing   again $a \ll_ {\mathbf  A} a'$, since \eqref{a1} and \eqref{a2} 
hold in $\mathbf A$ and $\ll_ {\mathbf  A}$  is
transitive. 

The proof that $\mathbf  D$ satisfies \eqref{a1} and
 \eqref{a2}, too, presents no significant variation with respect to 
the proof in (I) that $\mathbf  D$ satisfies  \eqref{a2}.
If, say, $x \leq y$ and $y\ll z$, then these relations are witnessed
by chains of relations either of the form  $u \ll_ {\mathbf  X} v$
or $u \leq_ {\mathbf  Y} v$, for $\mathbf X$ and $\mathbf Y$ 
either equal to $\mathbf A$ or to $\mathbf  B$, and at least one
 $\ll_ {\mathbf  X} $ appearing among the relations witnessing $y\ll z$.
Then a reduction similar to the one used in the proof of \eqref{a2}
in (I) shows that we can reduce ourselves to at most two relations,
with  $\ll_ {\mathbf  X} $ appearing at least once.
Indeed, adjacent $\mathbf A$-subscripted relations
can be evaluated in $\mathbf A$, hence can be reduced to a single
relation, using either \eqref{a1}, or \eqref{a2},
or transitivity of   $\leq_ {\mathbf  A}$ or of $\ll_ {\mathbf  A}$. 
Similarly for adjacent $\mathbf B$-subscripted relations.
If we have a chain of the form
\begin{equation}\labbel{red}    
 w \trianglelefteq _ {\mathbf  A} v \blacktriangleleft_ {\mathbf  B} u
\sqsubseteq _ {\mathbf  A} t,
  \end{equation} 
 then necessarily $v,u \in A \cap B =C$,
hence  $v \blacktriangleleft u$ can be evaluated in $\mathbf A$,
hence  \eqref{red} reduces to a single relation in $\mathbf A$. 
Since in the above reductions we use \eqref{a1}, \eqref{a2}
and transitivity of $\ll$, then at least one $\ll$ appears
in the reduced expression, since, as we mentioned, at least
one $\ll$ appears in the original relations witnessing $y\ll z$.

The proof in the above paragraph also gives
transitivity of $\ll$.
We now 
comment the difference between the present case and case (I).
 In case (I) we used (c1) in order to prove 
transitivity of $\ll$;  
we now explain why the additional 
assumption (c1) is not necessary when we have both
 \eqref{a1} and
 \eqref{a2} at our disposal.
 The typical example is given by 
 $x,z  \in A$, $y \in B$, with $x \ll y  \ll z$ 
given by $ x \ll_ {\mathbf  A} c \leq_ {\mathbf  B} y$ and
$ y \ll_ {\mathbf  B} d \leq_ {\mathbf  A} z$.
Notice that if we want \eqref{a2} to hold in $\mathbf  D$,
the above relations imply  $x \ll  y$ and $ y  \ll z$;
in other words, the satisfaction of \eqref{a2}
necessarily requires that  ${\ll_ {\mathbf  A}} \circ {\leq_ {\mathbf  B}}$
and 
${\ll_ {\mathbf  B}} \circ {\leq_ {\mathbf  A}}$ 
 appear in the definition of $\ll$ in $\mathbf  D$. 
On the other hand, from the above conditions we get
$  c \leq_ {\mathbf  B}  y \ll_ {\mathbf  B} d $,
and we need \eqref{a1} (not \eqref{a2}!) in $\mathbf  B$ 
in order to get  $  c \ll_ {\mathbf  B} d $, hence, as usual,
$  c \ll_ {\mathbf  A} d $.
Now from 
$ x \ll_ {\mathbf  A} c \ll_ {\mathbf  A} d \leq_ {\mathbf  A} z$
we get $ x \ll z$ by transitivity of $\ll$ and \eqref{a2} in $\mathbf A$.
However, we needed \eqref{a1} in order   to perform the above
computations, hence the argument cannot be carried over in case (I).

Now we deal with antisymmetry.
 Suppose that 
$\ll$ is antisymmetric in $\mathbf A$ and $\mathbf  B$ 
and $x \ll y \ll x $.
If both $x$ and  $y$ are in $A$, then 
$x=y$ by antisymmetry of $\ll$ in $\mathbf A$,
and similarly for $\mathbf  B$.   
Hence suppose, say,
$x \in A $ and $y \in B $. 
We first deal with the case when $x \ll y \ll x $
is witnessed by
$ x \ll_ {\mathbf  A} c \ll_ {\mathbf  B} y$
and 
$ y \ll_ {\mathbf  B} d \ll_ {\mathbf  A} x$, for $c,d \in C$,  
then we will show that all the other cases can be reduced to the
preceding case.
From 
$c \ll_ {\mathbf  B} y \ll_ {\mathbf  B} d$
we get  $c  \ll_ {\mathbf  B} d$
by transitivity of $ \ll_ {\mathbf  B} $
and then  $c  \ll_ {\mathbf  C} d$.
Symmetrically 
$d  \ll_ {\mathbf  C} c$,
hence $c=d$, by antisymmetry of $\ll_ {\mathbf  C} $.
Thus $ x \ll_ {\mathbf  A} c=  d \ll_ {\mathbf  A} x$,
hence $x=c$, by antisymmetry of $\ll_ {\mathbf  A}$.
Symmetrically, $y=c$, thus $x=y.$    
Now suppose that $x \leq_ {\mathbf  A} c$
appears among the relations witnessing  
 $x \ll y $.
Since we have already proved that $\ll$ is transitive,
then $x\ll x$, and this relation can be evaluated in $\mathbf A$,  
since we have already showed that $\mathbf A$ embeds in $\mathbf  D$.
From $x \ll_ {\mathbf  A} x \leq_ {\mathbf  A} c$
we get $x \ll_ {\mathbf  A} c$ by \eqref{a2} in $\mathbf A$.
Symmetrically, if $d \leq_ {\mathbf  A} x$,
then $d \leq_ {\mathbf  A} x \ll_ {\mathbf  A} x $,
hence $d \ll_ {\mathbf  A} x$, by \eqref{a1}. Performing the symmetrical
computations with respect to $y$, we can reduce all the $\leq$ 
relations to $\ll$ relations. This means that we always can reduce
ourselves to the case we have already treated.    

If $\ll$ is finer (coarser) than $\leq$ in all the models to be
amalgamated, then, by the definitions,
$\ll$ is finer  (coarser) than $\leq$ in $\mathbf  D$. 
We have completed
the discussion
of all cases.

The last two sentences in the theorem 
follow from \cite[Theorems 7.1.2 and 7.4.1]{H}
and \cite[Fact 2.1(3)]{KS}.
Compare the final part of the proof of Theorem \ref{cp}(1)(2)
at the end of Section \ref{2}.  
\end{proof}

\begin{remark} \labbel{add}
(a) The amalgamation property in Theorem \ref{gen}
still holds if we add families of $\leq$-preserving  
unary operations, and families of both $\leq$- and
$\ll$-preserving  
unary operations.

Indeed, the fact that $  D$ can be taken over the union
of (the images of) $ A$ and $  B$ implies that any two
unary operations $f_ {\mathbf  A}$ and $f_ {\mathbf  B}$
  defined, respectively,  on $\mathbf A$ and $\mathbf  B$ (and
agreeing on $\mathbf  C$)
can be extended uniquely to an 
operation $f$ on $\mathbf  D$.
The fact that we have strong amalgamation
warrants that $f$ is well-defined. 

Furthermore, superamalgamation implies that 
$f$ is $\leq$-preserving on $\mathbf  D$, since if, say,
$a \in A$, $b \in B$  and $a \leq_ {\mathbf  D} b$,
then  $a \leq_ {\mathbf  A} c \leq_ {\mathbf  B}  b$,
for some $c \in C$, thus
 $ f_ {\mathbf  A}(a) \leq_ {\mathbf  A} f_ {\mathbf  A}(c) $
and $ f_ {\mathbf  B}(c) \leq_ {\mathbf  B}  f_ {\mathbf  B}(b)$,
which imply
$f(a) \leq f(c)$ in $\mathbf  D$, since $c \in C$,
hence $f_ {\mathbf  A}(c) = f_ {\mathbf  B}(c)$. 
Similarly, the definitions of $\ll$ on $\mathbf  D$ 
in the proof of Theorem \ref{gen} imply that 
if    $f_ {\mathbf  A}$ and $f_ {\mathbf  B}$ are 
both $\leq$- and
$\ll$-preserving, then their  unique common extension on $\mathbf  D$
is $\leq$- and
$\ll$-preserving.

Many similar arguments are presented in \cite{duerel,apu}.
 
(b) The above remark does not hold for 
functions which preserve only 
$\ll$.
Let $\mathbf  C$ be a $3$-element poset
with three $\leq$-incomparable elements
$c$, $c_1$ and  $c_2$ and let no $\ll$-relation holds. 
Extend $\mathbf  C$ to $\mathbf A$ by   
adding two new elements
$a_1$ and  $a_2$
with $c_i \leq a_i$,  $c_i \ll a_i$, $i=1,2$
and no other nontrivial relation.
 Extend $\mathbf  C$ to $\mathbf B$ by   
adding two new elements
$b_1$ and  $b_2$
with 
$b_1 \leq c_1$, $c \leq b_2$
and  no $\ll$-relation.
Thus $\mathbf A$, $\mathbf  B$ and $\mathbf  C$ 
are posets with an auxiliary relation, actually,
causal sets (case $P = \{2,5 \}$, $Q = \{4 \}$, $N = \{ A1, A2, F \}$).

Add to each
model above a unary function
$f$ defined by $f(c)=c$, $f(x_1)=x_2$
and   $f(x_2)=x_1$, for $x_i=a_i,b_i,c_i$,
thus $f$ is $\ll$-preserving in 
all the above models, but not 
$\leq$-preserving in $\mathbf  B$.  
In any amalgamating model in which \eqref{a1}  holds, we should have
$ \kappa  (b_1)\ll \iota(a_1)$, hence,
if we want   $f$ to be $\ll$-preserving,
also $\kappa  (b_2) = \kappa  (f(b_1)) =   f(\kappa  (b_1))  \ll 
  f(\iota (a_1))   = \iota (f(a_1)) = \iota (a_2)  $.
Since $\kappa$  is supposed to be an embedding,
we should have $ \kappa (c) \leq \kappa (b_2)$, as well. 
By  \eqref{a1}, and since $\iota$ and $\kappa$ should agree on $C$,
 $ \iota (c) =  \kappa (c) \ll \iota (a_2)$.
 But $c \ll a_2 $ does not hold in $\mathbf A$,
hence $\iota$ fails to be an embedding. 
\end{remark}   

 Some assumptions on $P$, $Q$ and $N$ are necessary in
Theorem \ref{gen}. See Remark \ref{nec} below.

\section{Ordered structures with the superamalgamation property (continued)}
 \labbel{4}

In this section we show that a version of Theorem \ref{gen}
applies, more generally, to classes of ordered structures,
provided the superamalgamation property holds.
We will always assume that $\mathcal K$ is a nonempty class of
ordered structures of the same type, and $\mathcal K$
is closed under isomorphism.
The proofs in the previous section  are heavily  used here. 

Recall the conditions \eqref{a1}, \eqref{a2}
and the Conventions \ref{conv} from the previous section.  

\begin{lemma} \labbel{lem}
Suppose that $ Q \subseteq \{ 2,4,5 \} $ and
 $N \subseteq \{F, C, A1, A2 \} $.

Let $(D, {\leq_D)}$ and $(E, {\leq_E})$
be posets, and suppose that   $(E, {\leq_E})$ extends
 $(D, \allowbreak {\leq_D})$.
If $\ll_D$ is a binary transitive relation on $D$ such that    
$Q$ and $N$ are satisfied in $(D, {\leq_D},{\ll_D})$, then 
there is a transitive relation $\ll_E$ on $E$ such that
$Q$ and $N$ are satisfied in $(E, {\leq_E},{\ll_E} )$ and moreover
 $(E, {\leq_E},{\ll_E} )$ extends $(D, {\leq_D},{\ll_D})$.
 \end{lemma} 

\begin{proof} 
Again, we will consider  various cases.
 
(0) $A1 \notin N$ and $A2 \notin N$.
 This case is divided into two subcases. (0i) $C \notin N$. The subcase is
elementary;
it is enough to let $e \ll_ { E} f$
on $ E$ if and only if 
$e,f \in D$ and $e \ll_D f$
in $ D$
(if $2 \in Q$, in addition we let $e \ll_ {  E} e$, for all
$e \in E$). By construction, 
$\ll_E$ extends $\ll_D$.
Since $\ll_D$ is transitive on $  D$, then the new 
relation    $\ll_ { E}$ is transitive.
Similarly, if $\ll_D$ is 
antisymmetric, then $\ll_E$ is antisymmetric.
Suppose that
 $F \in N$,
that is, $\ll$ is required to be finer that $\leq$.
If this holds in $\mathbf  D$, then it holds in
$\mathbf E$, as well.

(0ii) $C \in N$. 
This is essentially the content of Lemma \ref{lemsolo}.
Just notice that antisymmetry of $\ll_ { D}$
is used there only in the proof of antisymmetry of $\ll_ { E}$.  
 
(I) $A1 \notin N$ and $A2 \in N$.
First, notice that if $C \in N$, that is, 
$ \ll $ is assumed to be coarser than $\leq $,
then \eqref{a2} (and also \eqref{a1}) are automatically
verified, if $ \ll $ is transitive.
Hence, in view of (0ii) above,
 in what follows we may assume that $C \notin N$.

Let 
$e \ll_ { E} f$
on $\mathbf E$ if  
$e \in D$ and there is $g \in D$
such that  $e \ll_ { D} g \leq_ { E} f$
(again, if $2 \in Q$, we 
additionally let $e \ll_ {  E} e$, for all
$e \in E$).
Condition \eqref{a2} is verified, since $\leq_ {  E}$ is
transitive. We have that $\mathbf E$ actually extends
$\mathbf  D$, since if $e \ll_D f$, then we get 
 $e \ll_E f$ by taking $g=f$. In the other direction,
if $f \in D$ and $e \ll_E f$, as witnessed by   
$e \ll_ { D} g \leq_ { E} f$, then $ g \leq_ { D} f$,
since $(E, {\leq_E})$ extends
 $(D,  {\leq_D})$. Thus 
$e \ll_ { D} g \leq_ { D} f$, hence 
$e \ll_ { D} f$ by \eqref{a2} in $\mathbf  D$.  
In particular, if $\ll_ { D} $ is antireflexive, then
$\ll_ { E} $ is antireflexive, since then $2 \notin Q$ and
if  $e \ll_ { E} e$, then $e \in D$.
If $\ll_ { D} $ is finer than $\leq_ {D} $,
then $\ll_ {E} $ is finer than $\leq_ {E} $,
by transitivity of $\leq_ {E} $.

We now check transitivity of $\ll_ { E} $.
If $e \ll_ { D} g \leq_ {  E} f
\ll_ { D} h \leq_ {  E} k$, then
$e,g,f,h \in D$, hence we get $e \ll_ {  D} h $
by \eqref{a2} and transitivity of $\ll$ in   $\mathbf  D$, thus
$e \ll_ {  D} h \leq_ { E} f$,
that is, $e \ll_ {  E} f$. This proves transitivity 
of $\ll_ { E}$.

It remains to deal with antisymmetry when $5 \in Q$.
If   $e \ll_ { E} f \ll_ {  E} e$
is witnessed by  $e \ll_ { D} g \leq_ {  E} f
\ll_ { D} h \leq_ {  E} e$, then
$e,g,f,h \in D$, hence 
$e \ll_ { D} g \leq_ {  D} f
\ll_ { D} h \leq_ {  D} e$,
since $\leq_ {  E}$ extends $\leq_ {  D}$, thus
$e \ll_ { D}  f
\ll_ { D}  e$
by \eqref{a2} in $\mathbf  D$.
Then $e=f$ by antisymmetry of 
 $\ll_ { D}$.

(II) The case $A1 \in N$ and $A2 \notin N$
is proved symmetrically.

(III) $A1 \in N$ and $A2 \in N$.
The case is  similar to (I), 
setting
$e \ll_ { E} f$
on $\mathbf E$ if  
 there are $h,g \in D$
such that  $e \leq_ { E} h \ll_ { D} g \leq_ { E} f$
(as usual, if $2 \in Q$, we 
additionally let $e \ll_ {  E} e$, for all
$e \in E$). We point out the main differences.

Transitivity of $\ll_ {  E} $ needs only a small variation.
If 
\begin{equation}\labbel{witt}
    e \leq_ { E} h \ll_ { D} k \leq_ { E} f \leq_ { E} p \ll_ { D} q \leq_ { E} g
   \end{equation} 
then $k \leq_ { E} p$, hence $k \leq_ { D} p$, since $k,p \in D$. 
From $h \ll_ { D} k \leq_ { D} p \ll_ { D} q$ we get
$h \ll_ { D} q$ by \eqref{a2} and transitivity of $\ll_ { D}$.
Hence $e \leq_ { E} h \ll_ { D} q \leq_ { E} g$, what we had to show.

A small detail needs to be worked out also in case 
 $4 \in P$. 
If, by contradiction, $e \ll_ { E}  e$, for some $e \in E$,
that is, 
$e \leq_ { E} h \ll_ { D} g \leq_ { E} e$,
for some $h,g \in D$,
then $g \leq_ { E} e  \leq_ { E} h $,
hence  $g \leq_ { E} h $,
$g \leq_ { D} h $, then
$h \ll_ { D} g \leq_ { D} h $,
thus $h \ll_ { D} h $,
by \eqref{a2},
contradicting antireflexivity of
$\ll_ { D}$, that is, contradicting
 $4 \in P$.  

The proof that $\ll_ { E}$
is antisymmetric if $\ll_ { D}$
is antisymmetric  is more involved.
Suppose that $e \ll_ { E} f \ll_ { E} e$ is witnessed by
\eqref{witt} with $g=e$. Then
$k,p \in D$, $k \leq_ { E} p$, hence 
 $k \leq_ { D} p$, thus we get
$h \ll_ { D}  p$
from \eqref{a2} 
and $h \ll_ { D} k \leq_ { D} p$.
Symmetrically, 
$p \ll_ { D}  h$, hence $h=p$,
by antisymmetry of $\ll_ { D} $,
hence also $h \ll_ { D}  p=h$.  
Similarly, $k=q$, using \eqref{a1}. From 
\eqref{witt}, recalling that presently $g=e$,
we have  $k=q \leq_ { E} g=e \leq_ { E} h$,
thus  $k  \leq_ { D} h$.
Then from $h \ll_ { D} h$ and \eqref{a1} 
we get $k  \ll_ { D} h$.
Since also $h \ll_ { D} k$, we get 
$h=k$ by antisymmetry of $ \ll_ { D}$.  
Then \eqref{witt} reads 
$e \leq_ { E} h = k \leq_ { E} f \leq_ { E} p = h = k= q \leq_ { E} g=e$,
thus we get $e=f$ by transitivity and
antisymmetry of $\leq_ { E}$.  
\end{proof}

\begin{theorem} \labbel{genn}
Suppose that $\mathcal K$ is a class
of ordered structures and
$\mathcal K$ has the superamalgamation property.
 Let $ Q \subseteq \{ 2,4,5 \} $
and $N \subseteq \{F, C, A1, A2 \} $. 
Suppose further that either
  \begin{enumerate}[(a)]   
 \item 
 $\{ A1, A2 \}  \subseteq N$, or 
\item
$A1 \notin N$, $A2 \notin N$, or
\item   
[(c1)] $F \in N$.
  \end{enumerate} 

Let $ \mathcal K _{Q,N} $ be the class of expansions of
models of $\mathcal K$ 
in a language with a further binary relation
$\ll$, with the conditions that
$\ll$ is transitive, $\ll$ 
satisfies the properties in $Q$ and that the properties
in $N$ are satisfied.

Then $ \mathcal K _{Q,N} $ has 
the superamalgamation property with respect to $\leq$.
 \end{theorem}

\begin{proof} 
The proof presents no essential difference with respect to 
the proof of superamalgamation in the proof of
Theorem \ref{cp} (1) given at the end of Section \ref{2}.  
Here just consider  members of $\mathcal K$, resp.,  $ \mathcal K _{Q,N} $ 
in place of models of $T^-$, resp., $T$  and use the more general
Theorem \ref{gen} and Lemma \ref{lem} in place of
Theorem \ref{due} and 
 Lemma \ref{lemsolo}.
Notice that, in the cases at hand, $P= \{ 2, 5 \} $,
hence clause (c2) in Theorem \ref{gen} is automatically satisfied.  
\end{proof}  

\begin{proof}[Proof of Theorem \ref{cp} (continued)]
Amalgamation for Part (3) is the special case
  $Q= \emptyset $  
 and $N= \{ F, A1,A2 \} $
 of Theorem \ref{genn}. For the theory of causal sets just take
$Q= \{ 4 \} $ instead. The proof at the end of Section \ref{2}
for the existence of Fra\"\i ss\'e limits, etc., works in the
present cases, as well. 
\end{proof}

\section{Urquhart doubly ordered sets}  \labbel{Usec}
  
According to
Urquhart \cite{U}, a doubly ordered set 
is a set $X$ endowed with two 
preorders $\leq_1$ and  $\leq_2$ on $X$ such that  
\begin{equation}\labbel{U}\tag{U} 
   \text{for all $x, y \in  X$, if 
$x \leq_1 y$ and  $ x \leq_2 y$,
 then $x = y$.}
   \end{equation}
Later, in the literature condition \eqref{U} 
is frequently omitted \cite{Hol}, henceforth
we will refer to an  \emph{Urquhart doubly preordered sets}
for a set with two preorders satisfying 
\eqref{U}. See the mentioned work \cite{Hol}
for more references and 
for a survey of applications of both the restricted
and the unrestricted notion.
Notice that we have explicitly mentioned
``preordered'' in the above definition; when we speak of a 
\emph{doubly ordered set} $X$ we will mean that 
$X$ is endowed with two partial orders (not just preorders).

\begin{theorem} \labbel{Uno}
The class of Urquhart doubly preordered sets
does  not have  the amalgamation property.
 \end{theorem}  

 \begin{proof} 
Let $\mathbf  C$ be a model with two elements
$c$ and $d$, both  $\leq_1$- and  $\leq_2$-unrelated.
Let $\mathbf A$ be a model 
with base set $A=C \cup \{ a \} $, where $a$ is a new element
distinct from $c$ and $d$, moreover, we assume that 
$a \leq_1 c$ and  $ a \leq_2 d$  hold in $\mathbf A$ 
and  no other nontrivial $\leq_1$ or  $\leq_2$ relation holds.
Let $B=C \cup \{ b \} $, where $b$ is still a new element,
$c \leq_1 b$ and  $ d \leq_2 b$  hold in $\mathbf B$ 
and, again,  no other nontrivial $\leq_1$ or  $\leq_2$ relation holds.
$\mathbf A$, $\mathbf  B$ and  $\mathbf  C$
are Urquhart doubly preordered sets, actually,
ordered sets.
 Let $\iota_1: \mathbf  C \to \mathbf A$ and
$ \kappa _1: \mathbf  C \to \mathbf B$ be the inclusions.

In any amalgamating model we must have 
$\iota(a) \leq_1 \iota(c) = \kappa (c) \leq_1 \kappa (b)$ and 
$\iota(a) \leq_2 \iota(d) = \kappa (d) \leq_2 \kappa (b)$,
since $\iota$ and $\kappa$ are supposed to be
embeddings, in particular, homomorphisms, and since
$\iota(c) $, $\iota(d)$ are supposed to be equal to,
respectively, $ \kappa (c)$ and 
$\kappa (d)$.
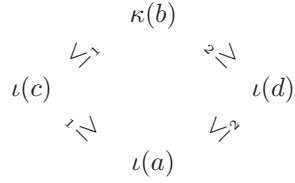
\begin{equation*}
\begin{array}{lcr} 
  & \kappa (b) & \cr
\enspace   \quad \quad  $\rotatebox[origin=c]{45}{$\leq_{_1}$}$  & & 
$\rotatebox[origin=c]{315}{${_{_2}}{\geq}$}$    \enspace  \quad  \quad  \cr 
\enspace   \iota(c)  & &  \iota(d) \enspace   \cr
\enspace    \quad  \quad  $\rotatebox[origin=c]{315}{${_{_1}}{\geq}$}$  & &
$\rotatebox[origin=c]{45}{$\leq_{_2}$}$  \quad   \enspace  \quad  \cr
 & \iota(a)
 \end{array} 
  \end{equation*}    

Hence, in any amalgamating model, $\iota(a)  \leq_1 \kappa (b)$ and 
$\iota(a) \leq_2  \kappa (b)$,
by transitivity of $\leq_1$ and $\leq_2$.
By \eqref{U},    $\iota(a)  = \kappa (b)$,
hence  $\kappa (b) = \iota(a)  \leq_1 \iota(c) = \kappa (c) $,
since $a  \leq_1 c$ and $\iota$ is supposed to be  
an embedding, in particular, a homomorphism.
Since $\kappa$ is supposed to be an embedding, from
$\kappa (b)    \leq_1 \kappa (c) $ we get
$b   \leq_1 c $, an inequality which fails in $\mathbf  B$,
thus there is no amalgamating model.
\end{proof}

In fact, the above proof shows more.

\begin{corollary} \labbel{corU}
The class of finite Urquhart  doubly ordered sets
does not have the amalgamation property
in the class of   sets with two transitive relations 
satisfying the Urquhart condition \eqref{U}. 

In particular, the class of Urquhart doubly ordered sets
does not have the amalgamation property.
In particular, the class of sets with two transitive relations 
satisfying the Urquhart condition \eqref{U}
does not have the amalgamation property.
 \end{corollary}

\begin{remark} \labbel{trans}
The use of transitivity is necessary in the proof of Theorem \ref{Uno}.
If $\leq_1$ and $\leq_2$ are assumed to be 
antisymmetric (or symmetric, reflexive, antireflexive),
but not necessarily transitive,
then, assuming, without loss of generality,
that the models to be amalgamated are 
$\mathbf  C \subseteq \mathbf A, \mathbf  B$   
with $A \cap B =C$, then a strong amalgamating structure
can be constructed over $D=A \cup B$ by setting
$\leq_i^ \mathbf  D = \leq_i^ \mathbf  A \cup \leq_i^ \mathbf  B$,
for $i=1,2$. This means that $x \leq_i^ \mathbf  D y$  if and only if
either $x,y \in A$ and $x \leq_i^ \mathbf  A y$,
or  $x,y \in B$ and $x \leq_i^ \mathbf  B y$.
In general, when, for relational languages,
a similar construction provides amalgamation, the class of models is 
frequently said
to have \emph{free amalgamation}. See, e.~g.,   \cite{F2};
however, notice that the terminology is not uniform in the literature.
When a class $\mathcal K$ has free amalgamation and, for some pair of 
binary relations, we consider the subclass of models 
in $\mathcal K$ satisfying
\eqref{U}, then free amalgamation is preserved, since 
if $x \leq_1 y$ and  $ x \leq_2 y$ in $\mathbf  D$,
then either $x,y \in A$, or $x,y \in B$.
In both cases, $x=y$, since $\mathbf A$ and $\mathbf  B$ 
are assumed to satisfy \eqref{U}.    
 \end{remark}

\section{Further remarks} \labbel{fur}

\begin{remark} \labbel{rmk}
In this remark  we show that
the superamalgamation property is necessary in 
Theorem \ref{genn}, case (b), even when $\ll$ is assumed to be
a partial order.
In other words, we do need Theorem \ref{josup} 
in the proof of Theorem \ref{cp}. 

(a) We first observe that
the theory $T_{uc}$ of partial orders which are 
expressible as  unions
of  chains has SAP.  
The theory $T_{uc}$ is a first-order universal theory:
just add to the theory of partial orders the sentence 
\begin{equation}\labbel{uc}      
\forall xyz (x \geq z \,\&\, y \geq z \Rightarrow x \geq y \text{ or } y \geq x),
 \end{equation}
 together with its dual

SAP for $T_{uc}$  follows from
SAP for the theory of linearly ordered sets.
For every chain in 
the amalgamating base $\mathbf  C$, 
amalgamate each corresponding triple of chains
as in the case of   linearly ordered sets.
If $\mathbf A$ and $\mathbf  B$ have more chains 
not connected to elements of 
 $\mathbf  C$,
just take their disjoint unions in 
the amalgamating model $\mathbf  D$.

(b) On the other hand, if we add another binary relation $\ll$
and add to $T_{uc}$  axioms asserting that $\ll$ is 
transitive and coarser than
$\leq$, then AP fails.

Suppose that $C= \{ c,d,e \} $, with all the
elements pairwise incomparable both with respect to 
$\leq$ and to $\ll$, and both $\leq$, $\ll$ reflexive. 
Let $A=C \cup \{ a \} $ 
with 
$d \leq a$, $d \ll a$, $c \ll a$, $a \ll a$ and 
no further $\ll$-relation, in particular,  
not $e \ll a$.   
Let $B=C \cup \{ b \} $ 
with 
$d \leq b$, $d \ll b$, $e \ll b$, $b \ll b$ and not  $c \ll b$.

If $\mathbf  D$ is an amalgamating structure 
which is a $\leq$-union of chains, then 
$a$ and  $b$ should be $\leq$-comparable,
say, $a\leq b$. If $\ll$ is coarser than
$\leq$ on $\mathbf  D$, then $b\ll a$.
If $\ll$ is transitive, then  $e\ll a$ in $\mathbf  D$,
but this contradicts the requirement that
$\mathbf A$ embeds in $\mathbf  D$, since
$e\ll a$ fails in $\mathbf  A$.

This shows that 
the assumption that $\mathcal K$ has the 
superamalgamation property 
is necessary in Theorem \ref{genn}.

(c) More generally, the counterexample in 
(b) above shows that the class of posets which are union of 
at most 3 chains 
and with a further coarser partial order
does not have AP in the class  
of posets which are union of chains and
with a further coarser transitive binary relation.
 \end{remark}

\begin{remark} \labbel{rmkk}
Similar arguments can be used to provide a
simple counterexample to results in \cite{sapimpap}. 

(a) We first observe that, for every $n \geq 1$,
the theory $T^n_{uc}$ of partial orders which are the union
of at most $n$ chains has SAP.  
Again,  $T^n_{uc}$ is a first-order universal theory:
just add to  $T_{uc}$ a sentence  asserting that there is no antichain of cardinality
$n+1$ (that is, for every $n+1$ elements, at least two of them
are comparable).

SAP for $T^n_{uc}$  follows from
SAP and JEP for the theory of linearly ordered sets.
If $\mathbf  C$ has already $n$ (distinct, nonempty) chains,
then no more chains can be present in $\mathbf A$ and $\mathbf  B$.
Similarly to Remark \ref{rmk}(a) above,
  amalgamate each corresponding triple of chains
as in the case of   linearly ordered sets.
If $\mathbf A$ and $\mathbf  B$ have more chains than $\mathbf  C$,
then amalgamate as above those chains having a representative
in $\mathbf  C$. Use JEP to deal with the other chains in 
$\mathbf A$ and $\mathbf  B$ in such a way that
the number of chains in the amalgamating model never exceeds $n$. 

(b) If $n >1$,  we add to the language a further unary operation $f$ and
assert that $f$ is isotone, then this extension of  $T^n_{uc}$
does not have AP.

Just let $C= \{ c \} $, $A= \{ c, a \} $,
$B= \{c, b_1, \dots , b_{n-1} \} $
with all the elements pairwise incomparable in each model. 
Let $f(c)= c, f(b_1)=b_1, \dots , f(b_{n-1}) = b_{n-1}$
and $f(a)=c$.   In any amalgamating model
$a$ should be comparable with one of the $b_i$'s
(since otherwise we have $n+1$ chains) 
but then  $f(a)=c$ and $f(b_i)=b_i$
contradict isotony, since $c$ and $b_i$ are
incomparable.   
 \end{remark}

\begin{remark} \labbel{noantich}
In contrast with  Remark \ref{rmkk}(a) above,
if $n \geq 2$, then the theory of partial orders with 
no antichain of cardinality $ > n$
fails to  have AP.   

Indeed,  let $C= \{c, c_1, c_2 \} $
with $c_1, c_2$ incomparable
and $ c \geq c_1 $, $c \geq  c_2 $. 
Let $A= \{a, c, c_1, c_2 \} $,
with $a \geq c_1$ and $a$ incomparable 
with all the other elements of $C$. 
Let
$B= \{ c, c_1, c_2, b_1, \dots , b_{n-1} \} $,
with the $b_i$'s pairwise incomparable, 
$b_i \geq c_2$, for every $i=1, \dots , n-1$,
and each $b_i$ incomparable 
with all the other elements of $C$.
In any amalgamating model
the set $\{ a, c, b_1, \dots , b_{n-1}\}$
has  $n+1$ elements, hence 
cannot constitute an antichain, so that
 $a$ should be  comparable with some $b_i$. 
If, say, $a \geq b_i$, then we get
$a \geq c_2$ by transitivity,  
but this is a contradiction, since $a \geq c_2$
fails in $\mathbf A$.
 \end{remark}   

\begin{remark} \labbel{nec} 
 Some assumptions on $P$, $Q$ and $N$ are necessary in
Theorem \ref{gen}.
For example, we prove that the  amalgamation property 
fails in case $P = \emptyset $, $  Q = \{ 5 \} $
and $N = \{F,A2 \} $ (thus (c2) does not hold).
 
We first exemplify the argument
showing that the strong amalgamation property fails.
Let $\mathbf  C$ be a one-element model
with $c \leq c$ and $c \centernot\ll c$.
Extend $\mathbf  C$ to $\mathbf A_i$ ($i=1,2$) by adding
 new elements $a_i \in A_i$ with $a_i\ll a_i$, $a_i\ll c$
and all the elements   $\leq$-related.
Then in all models \eqref{a2} holds,
$\ll$ is finer than $\leq$ and $\ll$ is 
transitive and antisymmetric.
If $\mathbf  D$ is an amalgamating model
with embeddings $\lambda_1$ and $\lambda_2$, then
from  $a_1\ll c $ and $c  \leq a_2$
we get $ \lambda _1(a_1)\ll \lambda _2(a_2)$, by \eqref{a2} and, 
 symmetrically, $ \lambda _2(a_2)\ll \lambda _1(a_1)$,
thus $ \lambda _1(a_1) = \lambda _2(a_2)$ by antisymmetry of 
$\ll$. Hence the strong amalgamation property fails.

The argument can be modified in order to show the failure
of the amalgamation property.
Consider $\mathbf  C$ as above, with a further element
$d$ such that $d \leq c$, not $c \leq d$
and no further $\ll$ relation holds.
Let $a_1\ll d$, $a_1\leq d$ in $\mathbf A_1$ and
not $a_2\ll d$ in $\mathbf A_2$.
The above argument shows  that
$ \lambda _1(a_1) = \lambda _2(a_2)$,
for any pairs of embeddings given by a realization
of the amalgamation property.
But if the $\lambda_i$s are embeddings, 
we get both 
$ \lambda _1(a_1) \ll \lambda _1(d)= \lambda _2(d)$ and
$ \lambda _1(a_1) = \lambda _2(a_2)\centernot\ll \lambda_2 (d)$,
a contradiction.
\end{remark}

\end{document}